# Information Surfaces in Systems Biology and Applications to Engineering Sustainable Agriculture

Hesam Dashti[1], Alireza Siahpirani[2], James Driver[1], Amir Assadi[1]

[1]Department of Mathematics and [2] Electerical and Computer Engineering, University of Wisconsin, USA
{Dashti, Fotuhisiahpi, Driver, Ahassadi}@wisc.edu

**Abstract.** Systems biology of plants offers myriad opportunities and many challenges in modeling. A number of technical challenges stem from paucity of computational methods for discovery of the most fundamental properties of complex dynamical systems. In systems engineering, eigen-mode analysis have proved to be a powerful approach. Following this philosophy, we introduce a new theory that has the benefits of eigen-mode analysis, while it allows investigation of complex dynamics prior to estimation of optimal scales and resolutions. Information Surfaces organizes the many intricate relationships among "eigen-modes" of gene networks at multiple scales and via an adaptable multi-resolution analytic approach that permits discovery of the appropriate scale and resolution for discovery of functions of genes in the model plant Arabidopsis. Applications are many, and some pertain developments of crops that sustainable agriculture requires.

## 1 Introduction

The concept of dynamical systems has been proposed to investigate natural complex time-dependent systems. Poincare first introduced dynamical systems to represent the topological features of orbits[1], and this notion of dynamical systems has been extended to numerical methods to model evolutionary systems that may, or may not, have cyclic behavior. As illustrated in [1], Poincare incepted a rich concept that allowed for the extension of the field of dynamical systems to capture many complex problems in computational systems biology [2-4], numerical analysis[5-8], number theory [9],and to classical and quantum mechanical [10, 11] complex systems.

A common method to analyze dynamical systems is to consider the system as a set of evolutionary functions in a geometrical manifold. However, this requires complex computations. These complexities have motivated researchers to adapt current methods, or develop new computational tools to facilitate these computations. Among these, soft computing [6], probabilistic and stochastic methods [3, 12, 13], topological and geometrical methods [14, 15], and machine learning [2, 8, 16-18] are most notable, and have brought a new level of understanding to the study of dynamical systems.

The numerical study of dynamical systems is focused on modeling the current state of the system [5, 6, 19] for data mining purposes (i.e. supervised and unsupervised classification) and proposing a model for future states through the introduction of predictive models [13, 18, 20]. Hence, dynamical systems may be considered as elements of a higher dimensional space. Explicitly, instead of investigating the finer features contained within a dynamical system, one can develop a model to explain the relationships amongst systems[7, 21]. Thus, introducing a measurement for quantifying the distance of two dynamical systems is necessary. Lack of such a measurement might be more sensible in the field of biological dynamical systems, because of the existence of the phenotypic and genotypic variation within species. However, measuring similarities or dissimilarities of biological dynamical systems is a mathematical representation of what is known as quantifying phenotypic traits of dynamical genotypic perturbations [22, 23].

One such biological dynamical system of interest is a time-series of gene expression profiles, where the data set is comprised of a set of genes stored in rows along time-steps corresponding to expression values in columns. The analysis of such matrices have included the aforementioned techniques and the reader is encouraged to examine[24-29] for a comprehensive review. Among these, Bar-Joseph Ziv et. al [28] introduced a novel method to classify dynamical systems through the introduction of a model that considers a discrete biological dynamical system as a continuous system to allow for the clustering of its current state to predict future states. Herein, we consider a dynamical system as a continuous system and introduce a novel methodology, InfoSurf (Informative Surfaces), to represent the system. Based on this representation, we propose a measurement for quantifying differences between dynamical systems. Moreover, InfoSurf is able to identify the objects responsible for the differences between the dynamical systems under consideration.

InfoSurf is designed with respect to three well-known mathematical theories, namely, multiscale analysis [30], multiresolution analysis [31], and Eigenanalysis [32] on an instantiation of the Sobolev space [33]. These techniques are applied on a high-performance computing platform where we consider a dynamical system as a two-dimensional table consisting of '$m$' objects in the rows and '$n$' time points in the columns. InfoSurf uses the value of the $i^{th}$ object at the $j^{th}$ time point as the height ($z$-coordinate) of a point ($x=i$, $y=j$) in a three-dimensional Euclidean space. Considering objects in this manner allows for the construction of a surface which is representative of the entire system and is the starting point of our analysis.

The novelty of InfoSurf and its contribution to the field of studying dynamical systems is described in section 2. In section 3 computational steps of InfoSurf's methodology are illustrated. In section 4 we evaluate the methodology on a biological dynamical system acquired from the shade avoidance study [34] on *Arabidopsis Thaliana*.

## 2 Contribution to Value Creation

In view of critical significance in economics and international relations, sustainable agriculture is regarded as a fundamental research domain that needs transformative innovations beyond the current incremental successes. Molecular methods in biotechnology and agricultural engineering promise rapid breeding of new lines of crops that would sustain stress from global warming and other harsh climatic events. Success of molecular methods depends on breakthrough in molecular systems biology, and invention of new ways of understanding the complex dynamics formed by time-course data from genes, proteins and other biomolecules. The technical demand for development of new algorithms to surmount the present computational challenges requires re-examination of traditional methods that have proved successful in non-complex systems and their dynamics. In particular, researchers must address discovery of the necessary biological properties implicit in –omic data, and mine the abundance of dynamical features that could be observed only in appropriate scales and via optimal resolutions.

This research addresses some of the bottlenecks that are posed in providing effective applications of systems biology to sustainable agriculture. Thus, the applications of this research will contribute towards value creation and directly addressing critical scientific problems that face humankind today.

## 3 Methodology

One of the novelties of InfoSurf is that it provides a new representation for a dynamical system. InfoSurf allows every dynamical system ($M_{mxn}$) to be considered as a surface in three-dimensional Cartesian coordinates, where $m$ objects are stored in $m$ rows whose corresponding $n$ time-dependent attributes are stored in $n$ columns. Biological dynamical systems, especially gene expression time-series, consist of $m$ genes with $n$ expression values where $m \gg n$. This is due to the fact that performing biological experiments for long periods of time with many points to assess expression values is expensive and inefficient [35]. Therefore, the number of time points is much smaller than $m$ and many applications have been developed for analyzing gene expression time series with $n \leq 8$ [35-37]. Hence, in the time-series under consideration there are few time points, and the interpolation of these points to include a greater number of attributes for each object is reasonable [28, 29, 38].

InfoSurf uses an interpolation method for incorporating more time points, extending the number of the columns, of the matrix $M_{mxn}$. This method is a row-wise interpolation and its algorithm was chosen with respect to two constraints: the interpolation error and the regularity of its functions. The fitted function to every row passes through the points of the row with an error of (-4.7726e-12). The regulatory of the fitted functions is important to ascertain the smoothness of the corresponding surface (this yields raw-wise smoothness). InfoSurf also requires regularity between different objects (column-wise regularities). To achieve such regularity, InfoSurf sorts the objects based on three features: the area underneath a) the signal (a row), b) its first derivative, and c) its second derivative. After sorting the data based on the comparison of these features, a row that has the closest similarity to its neighbors in terms of value, speed, and concavity, may not yield a smooth surface, but it provides the best orientation of the dynamical system for our analysis without introducing error. After these preliminary steps, InfoSurf projects the dynamical system onto a smoother surface without information loss through the use of multiscale and multiresolution analysis in conjunction with singular value decomposition (SVD).

### 3.1 Multiscale Analysis

Multiscale analysis has become prominent in recent years due to advances in computational speed and an increasing number of problems that rely on disparate mathematics to describe phenomena at different spatial and temporal levels. Multiscale analysis [30, 39] provides a bridge between these levels and allows one to analyze phenomena that are interdependent and might otherwise fail to be properly described within the scope of a single model. As problems become more sensitive to their representations, it becomes impossible to ignore phenomena simply because it does not appear to have a known source at a particular temporal or spatial resolution.

This is particularly true of systems biology; where the dynamical quantum [10] properties of physics yield the atoms of chemistry whose bonds go on to form the biological molecules of interest in the system of the organism under study. Hence, multiscale analysis plays an important role for scientists wishing to analyze data and probe measurements at different levels to identify new phenomena that may otherwise go undetected.

This ability is especially important for the continued advance of genetics where the size of a particular genome may be very large. Indeed, genetic systems biology has posed many new problems for data analysis due to the massive amount of information that must be analyzed for a given experiment. Multiscale analysis is well suited to address this problem since it is easily formulated to analyze databases and able to render data in many different fashions to provide key insights much faster than an individual researcher. It is for this reason we find multiscale analysis appropriate and develop a model herein to apply to multiple dynamical systems under unique conditions and develop a new measurement to compare them.

To perform multiscale analysis on a dynamical system ($\mathbf{M}_{mxn}$), InfoSurf considers a sliding window, a sub-matrix $S$, of size $k \times k$, $2 \leq k \leq \min(m, n)$, of $\mathbf{M}_{mxn}$. The size of $S$ varies between the construction of different surfaces but remains invariant for the entire surface under consideration and for the comparison of two surfaces as will be described later in this section. The sub-matrix slides in two directions; the first sub-matrix is defined by $S = \mathbf{M}(1:k, 1:k)$ (left-top), and slides to right and down by one in every iteration. The following pseudo-code illustrates the process:

```
for i=1:m-k
    for j=1:n-k
        S = M(i:i+k, j:j+k)
        //Performing analysis on S
    end
end
```

One finds that this process projects the matrix $M_{mxn}$ to a super-matrix containing ($M$-$k$)x($N$-$k$) sub-matrices of dimension $K$x$K$.

Considering overlap in our sub-matrices is a fundamental aspect of InfoSurf. Overlap reveals the continuous influence of objects on other groups of objects and allows for the method to proceed continuously and reveal information between data points that would otherwise be unaccounted for.

Considering every point in different windows illustrates the effect of an object on another object/objects and is seen multiple times as the object remains in the sliding window. This amplifies the effect(s) of the object and allows for it to be observed in different sliding windows. This allows for easier identification of an object and increases the accuracy of the algorithm when analyzing a dynamical system.

### 3.2 Multiresolution Analysis

Multiresolution analysis[40, 41] allows for larger features of a system to be reduced to the relationships of its fine features. For an example, in a gene expression time-series it allows for the detection of groups of genes that are potentially up or down regulated with respect to one another when verified through relevant biological data. Through use of surfaces, one can observe patterns of gene activity and reduce the macroscopic picture to the action of the individual genes responsible. Considering subsets of genes through different resolutions increases the accuracy of the InfoSurf algorithm and reduces run-time complexity. The different resolutions of InfoSurf are characteristic of the sliding window ($S_{kxk}$) described in the previous section. Through use of this window, InfoSurf's detection capabilities are increased and it allows for the extraction of specific attributes of genes and the construction of their interrelationships.

Starting multiscale analysis at a larger scale, larger $k$ for the size of the sliding window, allows the algorithm to identify regions of differences of two dynamical systems. In the multiresolution process, the regions of differences are then considered with smaller values of $k$ and the process reveals more details of the surfaces of the regions at a finer scale. InfoSurf uses the multiresolution process to zoom into the regions with very fine sliding windows and identify the specific objects corresponding to the differences between the dynamical systems.

Multiscale analysis allows one to examine the hierarchical structure of the objects from coarse to fine interrelations. It provides the ability to capture relationships between groups of objects (coarse scale) and tune it to identify relation between the objects (fine scale) [42].

### 3.3 Eigen Analysis

Eigen analysis is a fundamental method of data analysis and the investigation of structural properties of datasets [43-47]. The use of Eigen analysis in the InfoSurf algorithm was inspired by the kinematics of surface deformation as described in [48]. Here we describe the concept and our reasons for Eigen analysis in InfoSurf via examples.

**Formulation of Kinematics of Dynamical Systems**

InfoSurf considers a dynamical system as a surface in a three- dimensional space. To compare two dynamical systems, InfoSurf sorts one of the two surfaces based on the sum of its three representatives for each of its objects (e.g. genes); the integral of the signal, the integral of its speed, and the integral of its concavity. The one chosen sorted surface is a baseline (control surface, **X**) to be compared to another surface under consideration (deformed surface, *x*). Figure 1 shows an example of these surfaces in a three-dimensional space.

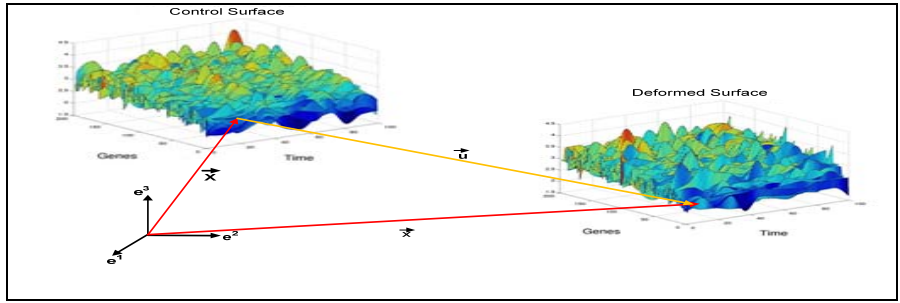

**Fig. 1.** Representing the deformation in dynamical systems. $e^1$, $e^2$, and $e^3$ are perpendicular basis vectors of a three-dimensional space. Every point on the control surface and the deformed space is represented by $\vec{X}$ and $\vec{x}$ in that basis space, respectively. The vector $\vec{u}$ shows the displacement of $\vec{X}$.

We define the vector $\vec{u}$ as the displacement of the vector $\vec{X}$. $x$ is the deformed surface and $\vec{x}$ is its position vector in the canonical basis ($e^1$, $e^2$, and $e^3$). This yields the following equation:

$$(1)\ \vec{u} = \vec{x} - \vec{X}$$

It is clear that position of a point in the deformed surface depends on $\vec{X}$ and the displacement $\vec{u}$, where every point in $X$ has its own displacement vector $\vec{u}$. Hence, the position of a point in x, is a function of the corresponding point $X$, $\vec{X}$, and its displacement $\vec{u}$. In other words, $\vec{x}(\vec{X}) = \vec{X} + \vec{u}(\vec{X})$. Since the control surface is continuous, one can consider infinitesimal changes in $\vec{X}$ and consider projections of small portions of the control surface into the deformed space, as shown in Figure 2.

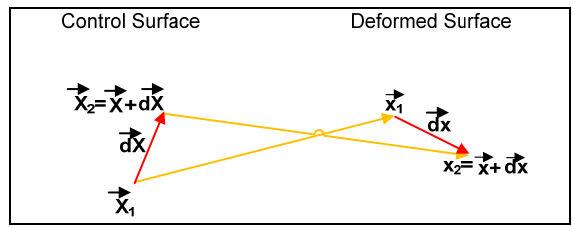

**Fig.2.** This graph shows the infinitesimal changes from the control surface into the deformed surface.

Since $d\vec{x}$ depends on $\overrightarrow{dX}$ and changes between $\overrightarrow{x_1}$ and $\overrightarrow{x_2}$, the relation between $\overrightarrow{dX}$ and $d\vec{x}$, when $\overrightarrow{dX} \to 0$, can be written as

$$(2)\ \overrightarrow{dx} = \overrightarrow{dX} \cdot \nabla \vec{x}$$

where $\nabla \vec{x}$ is the Jacobian matrix or the position gradient, defined as follows:

$$J = \nabla \vec{x} = \frac{\delta x}{\delta X}.$$

The projection of the Cartesian space into a spherical space and considering the vectors $\vec{x}$ and $\vec{X}$ in the spherical space, the Jacobian matrix can be calculated as:

$$J_{ij} = \frac{\overrightarrow{x_j}}{\overrightarrow{X_i}} = \frac{d\vec{x}}{dr_j}\frac{dr_i}{d\vec{X}}.$$

The InfoSurf algorithm uses the Jacobian matrix as a measurement of the dissimilarities between the dynamical systems. In [48] it has been shown that the sum of the eigenvalues of the Jacobian matrix corresponds to the changes between two surfaces, where if the sum is positive it means the amount of changes corresponds to less deformation, zero implies the surface has been crushed to a point, and a negative sum is representative of twists developing in the surface.

InfoSurf, instead of calculating this measurement for the entirety of the surfaces under consideration, uses multiscale and multiresolution analysis to capture the differences of many small covering sub-surfaces of the dynamical systems. For every invariant sliding window InfoSurf computes its eigenvalues and eigenvectors. Considering $Sv_i = \lambda_i v_i$, the eigenvalues, $\{\lambda_1, \lambda_2, \ldots, \lambda_n\}$, are the confidents of the corresponding eigenvectors $v_i$.

Where $S$ is the sliding sub-matrix of the data matrix $\mathbf{M}_{mxn}$. If we consider the matrix $S$ as a surface, where $S(i, j)$ represents the heights of the surface at $(i, j)$, the eigenvalues represent **a)** the heights of the surface. The distribution of eigenvalues is representative of the number of eigenvectors needed to reconstruct the $S$. **b)** Eigenvalues represent the smoothness of the surface. The following figures illustrate these properties.

Consider two matrices **A** and **B**; **A**$(i, j) = 1$, **A**$(i, i) = 2$ and **B**$(i, j) = 2$, **B**$(i, i) = 4$, where $1 \leq i, j \leq 100$ (Figure 3).

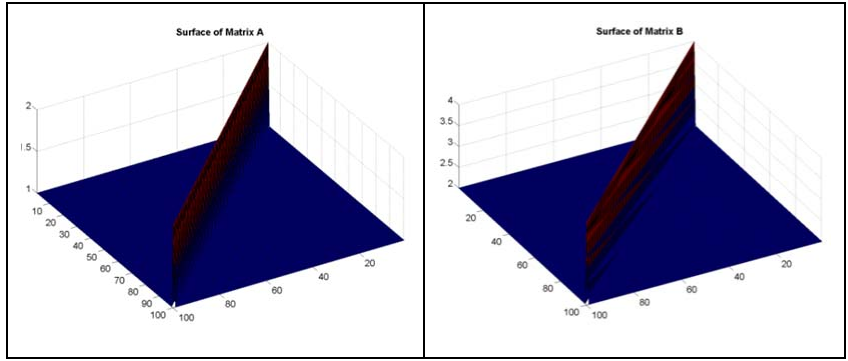

**Fig.3.** This figure shows the surface representation of A and B. The value of the main diagonal of matrix A is 2 and all other elements are 1. The value of the main diagonal of matrix B is 4 and all other elements are 2.

For every dynamical system InfoSurf computes the eigenvalues of every sliding window. Since the sliding window iterates in two-dimensions the eigenvalues are stored in a matrix, **E**. For the $r^{th}$ row and $s^{th}$ column iterations, InfoSurf associates the absolute value of the sum of the eigenvalues to **E**($r$, $s$). The matrix **E** is called an Eigensurface and represents the internal properties of data. Figure 4 shows the Eigensurfaces of the matrices **A** and **B**.

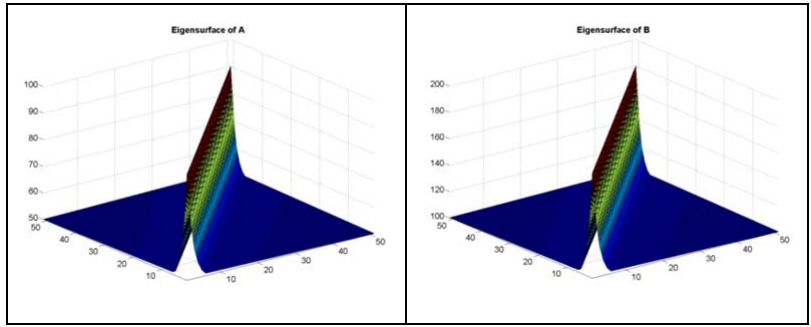

**Fig.4.** Eigensurfaces of A and B. The heights of these surfaces represent the heights of the corresponding surfaces in Fig 3. The shape of the two Eigensurfaces is the same but the heights are different. The Eigensurface of B is twice the height of the Eigensurface of A.

One can observe the differences in the heights as differences of the heights of the original matrices A and B. Normalizing the Eigensurfaces cancels the effects of the differences of heights. This is shown in figure 5.

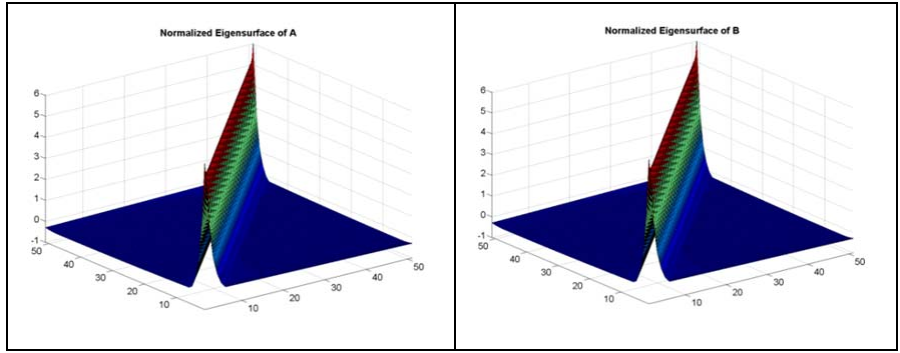

**Fig.5.** This figure shows the normalized Eigensurfaces of figure 4. The two surfaces are now identical.

In order to see how eigenanalysis captures smoothness of the surfaces, another example is given in figure 6, where two different matrices of the same height and with different structure are shown.

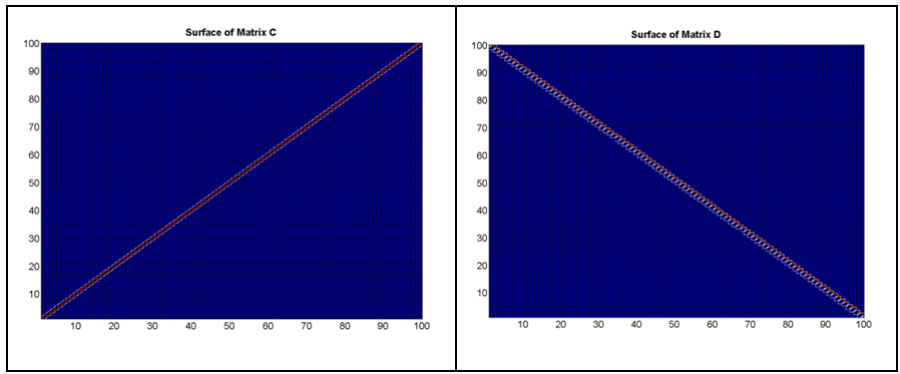

**Fig.6.** Surfaces of two matrices C and D are shown, where C(i, j) = 1, C(i, i) = 2, D(i, j) = 1, D(i, 100-i+1) = 2, $1 \leq i, j \leq 100$. This figure shows the surfaces of two matrices C and D. The main diagonal of matrix C is 2 and all other elements are 1, while the anti-diagonal of matrix D is 2 and all other elements are 1.

After constructing Eigensurfaces of C and D and normalizing them, one can observe the disparate surfaces as shown in figure 7.

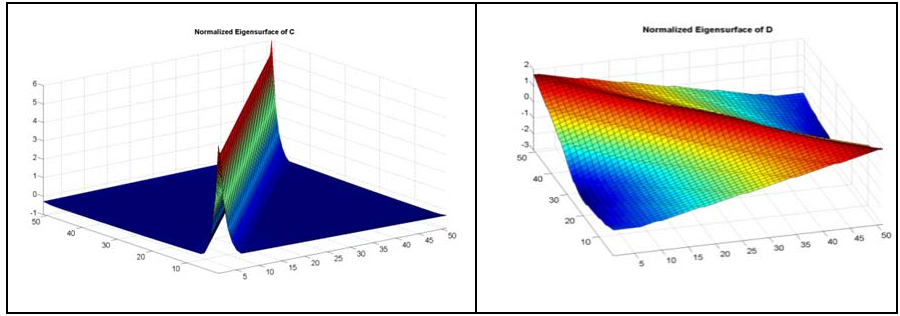

**Fig.7.** Eigensurfaces of matrices C and D after normalization. The height of the two surfaces (from minimum to maximum) is identical, but differences in the original surface caused the changes in the Eigensurface.

Recall that InfoSurf considers the first and second derivative of a dataset. After constructing the Eigensurface, InfoSurf calculates the first and second derivatives of the surface. These derivatives are useful for identifying circadian clock information of dynamical systems [49, 50]. The derivatives are calculated by a second order approximation method [51].

The first order mixed derivative of a surface allows for the association of a measure to the change of information content within a particular window of the objects under consideration. This describes the amount that an object within this window perturbs another object, or objects, within the window. While the first derivative is characteristic of the slope of the change of the eigenvalues, relating the change of the information content of each window and the objects within it, the second derivative provides information on the concavity, or acceleration of changes, circadian clock, and shows whether a subset of objects within each window is having a larger or smaller effect as time progresses. This allows for the identification of a subset of objects whose local maximum or minimum from their representative eigenvalues are characteristic of one of the objects within the window being responsible for the perturbation of the other objects (i.e. an increase or decrease in time-dependent attributes of other objects). Figure 8 shows the derivative surfaces of the Eigensurfaces of the matrices C and D as introduced in Figure 7.

The outcomes of analyzing the derivatives are independent of the ordering of the objects since multiresolution accounts for all window sizes, and despite row exchange, will carry the same perturbation. Its first and second derivatives will still retain their properties of detecting maxima and minima since the information is not lost but trans-located to another region of the matrix. Hence, another region of the matrix will contain the same maxima/minima peaks as before, simply in a different location of the surface.

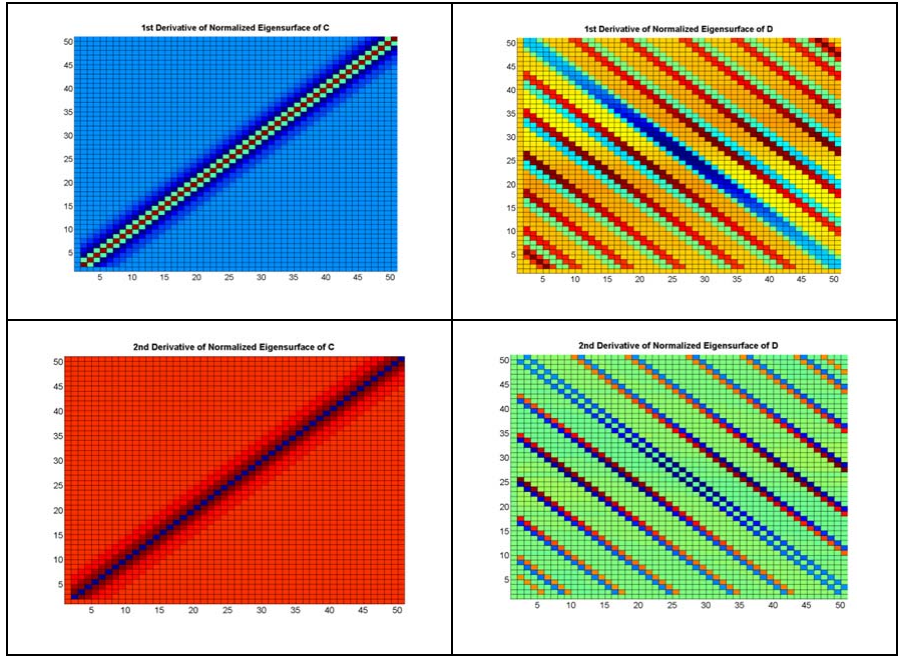

**Fig.8.** This figure shows the first and second derivatives of the two normalized Eigensurfaces. These surfaces illustrate how the Eigensurfaces are changing.

After constructing the representative surfaces, Eigensurfaces along with their first and second derivatives, and the Jacobian matrix, InfoSurf measures the dissimilarity between the dynamical systems. This dissimilarity can be represented by the Jacobian matrix (represented by a surface) and the distance and free-scale distance of the derivative surfaces.  Therefore to compare two dynamical systems, InfoSurf generates seven surfaces: a) distance of the Eigensurfaces, the surfaces of the first derivatives, and surfaces of the second derivatives.  b) free-scale distance of the three representative surfaces, and c) the Jacobian matrix.  The distance is the absolute value of the direct subtraction of two matrices (surfaces):

$$Dist(A, B) = abs(A - B)$$

and the scale-free distance is defined :

$$FreeDist(A, B) = \frac{abs(A - B)}{abs(A) + abs(B)}$$

The distance surfaces show differences/similarities of the dynamical systems. Figure 9 shows *Dist* and *FreeDist* of the derivative surfaces of Eigensurfaces of **C** and **D**. Flow of data analysis in the InfoSurf algorithm is depicted in figure 10.

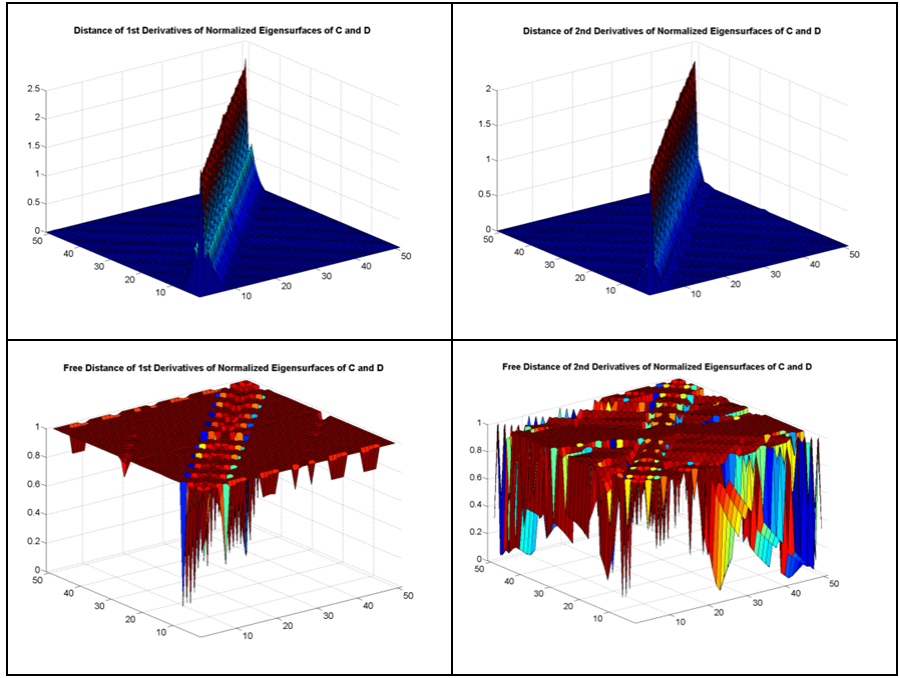

**Fig.9.** This figure represents the Distance and Free Distance of the first and second derivatives of the Eigensurfaces shown in figure 8. The Distance surface highlights the largest differences in the original surfaces. The difference surface is then divided by the sum of the heights of original surfaces to give the Free Distance. Free Distance highlights the surfaces' differences when compared to global maxima. While some changes appear small, they are significant locally.

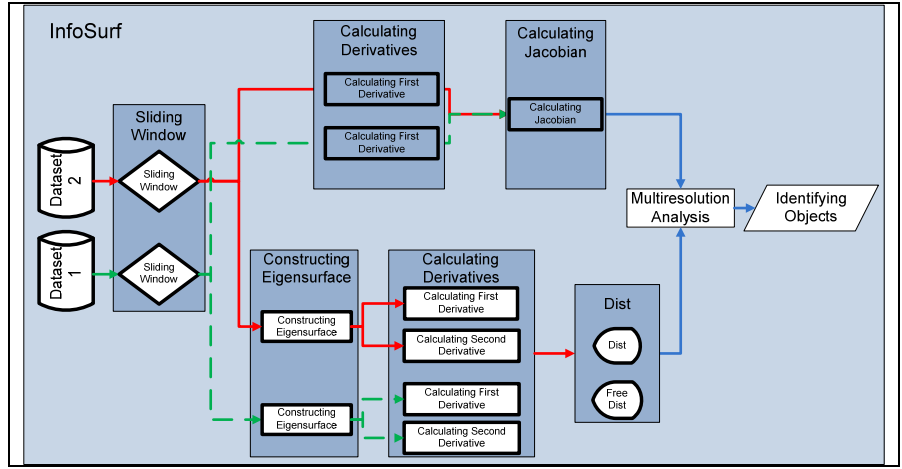

**Fig.10.** InfoSurf process. This figure shows the flow of data in the algorithm.

As we mentioned earlier, the derivative surfaces are useful for demonstrating the characteristics of behavior of a surface. After calculating the Distance, Free Distance surfaces, and finding the largest differences between the derivative surfaces, one can locate the biggest differences in behavior of Eigensurfaces and therefore approaches towards an area in the original data set that caused these differences. In order to identify the objects whose are responsible for the differences between the surfaces, InfoSurf performs a multiresolution analysis on the regions with largest differences and in this way narrow its analysis. An example of InfoSurf different steps is given in the next section.

## 4 Discussion and Experimental Results

To evaluate the InfoSurf method we used two dynamical systems from the shade avoidance experiment [34]. The data represents gene expression levels of *Arabidopsis Thaliana* when stimulated with changes in temperature and light. Samples of approximately 23 thousand genes were taken in 4 hour intervals over a period of 48 hours. We chose two specific data sets to compare. The first experiment consisted of exposing the plants to constant light and a rise in temperature for 12 hours and then 12 hours of darkness with a lower temperature. The second experiment was from plants that were exposed to light during the entire experiment while the change in temperature was the same as the previous experiment. The first data set is called LDHC (Light, Dark, Hot and Cold) and the second data set is called LLHC (Light, Light, Hot and Cold). Each data set has 22,810 rows (the number of studied genes of *Arabidopsis Thaliana*) and 12 columns (the number of four hour time steps in the experiments). We interpolated the data sets row-wise to obtain 100 points for each gene.

To acquire a smoother starting surface, we sorted the genes through row exchange based on the similarity of their time series (i.e. expression values). If we denote the time series of a gene by f(t), we consider the value of $g(x) = \int f(x) + \int f'(x) + \int f''(x)$ to be a good representation of the shape of the signal. The integrations have calculated by the trapezoidal approximation. By sorting the genes according to this value we obtain a smoother surface. However, if we sort both data sets with this method, the order of the genes might be different. Therefore we sort only one data set (control surface; LDHC) and rearrange the other to have the same order of genes as the smoothed one (deformed surface; LLHC). Figures 11 and 12 demonstrate the original and sorted data sets respectively.

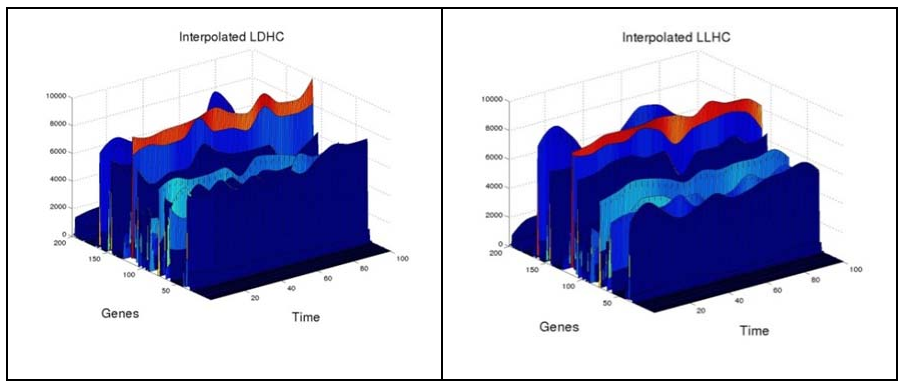

**Fig.11.** The first 200 genes of interpolated data sets. This figure shows the surfaces generated from the first two hundreds genes of the two data sets LDHC and LLHC. The original order of the genes, some genes with small expression values and some genes with larger expression values are adjacent.

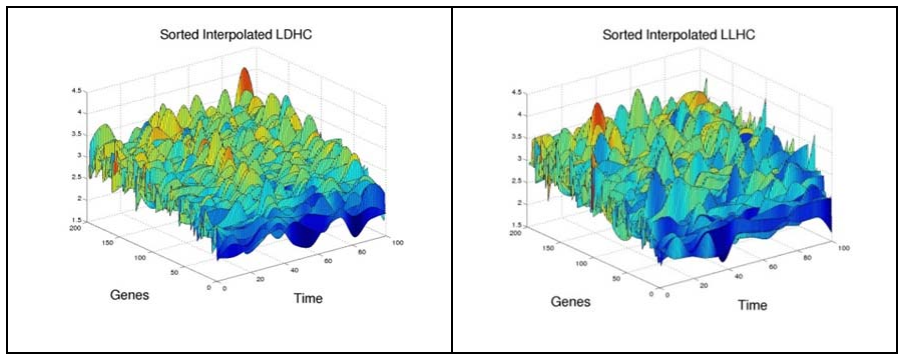

**Fig.12.** Sorted data sets of the first 200 genes. This figure shows the surfaces generated from two data sets after sorting one of the sets (deformed; LLHC) based on the smoothed one (control; LDHC). As illustrated, the heights of adjacent genes are in the same value range.

InfoSurf starts with multiscale analysis by a sliding window ($S_{40 \times 40}$), and calculates the Eigensurfaces for both data sets and their individual first and second derivatives. Figure 13 shows the Eigensurfaces of the sorted data sets, and figure 14 demonstrate first and second derivative of these surfaces.

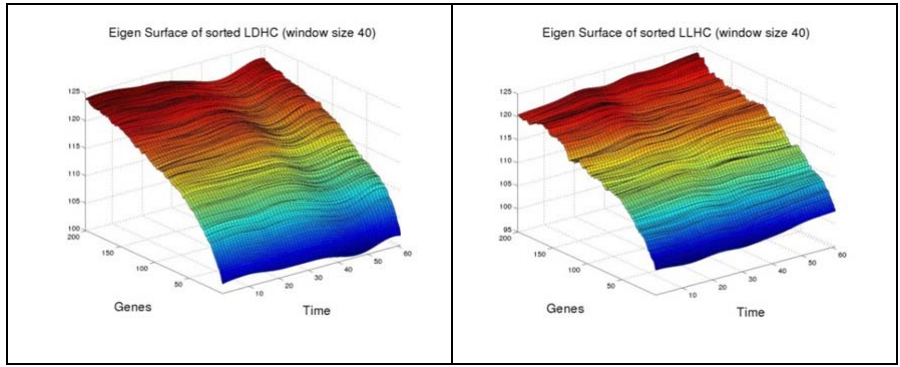

**Fig.13.** This figure shows the Eigensurfaces of the sorted data sets. The change in height of surfaces represents the change of heights in the original surfaces (shown in figure 12).

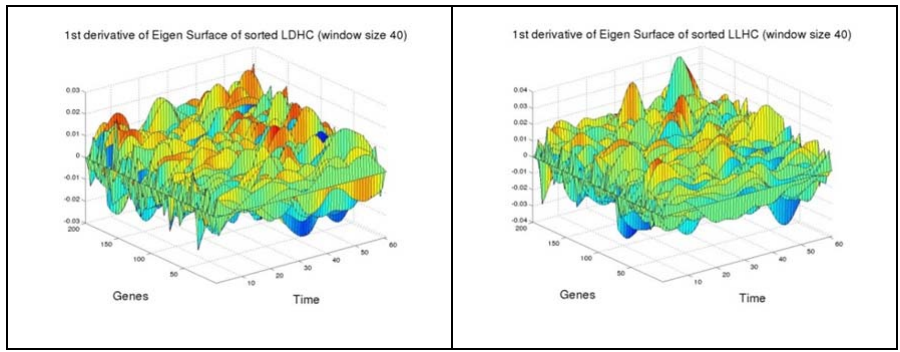

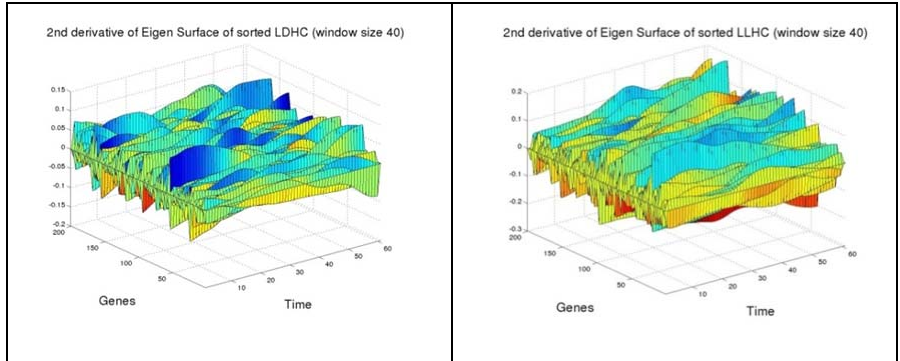

**Fig.14.** This figure shows the first and second derivative of two Eigensurfaces.

To locate the genes that act differently in the two data sets, we considered the differences of second derivative of the two Eigensurfaces with the sliding window of size 40 (figure 15 shows the Dist and FreeDist of those derivative surfaces) and found the local extrema. These points represent a window of 40x40 in the original data sets (40 genes in 40 time steps) whose eigenvalues are different between the two data sets. To refine the selection of genes, we used a higher resolution sliding window (20x20) to consider inside of the 40x40 matrix. The Eigensurface is constructed and the second derivative is calculated have a better understanding of the genes behavior. This leads to a 20x20 window in the original data that includes the local exterma. We resumed the increase of resolution by using the Eigensurfaces of the10x10 and 5x5 sliding windows which yield a 5x5 area (5 genes in 5 time steps).

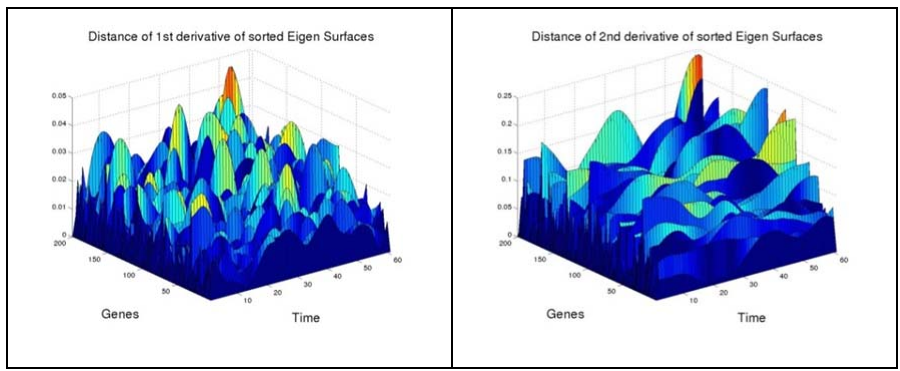

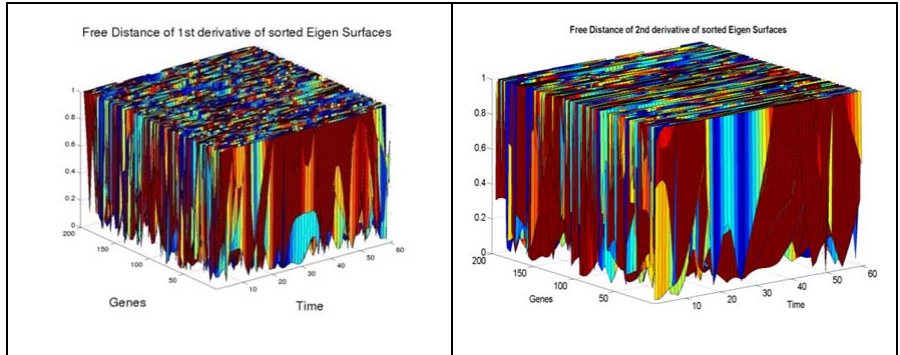

**Fig.15.** This figure shows the absolute value of the Distance and Free Distance of the first and second derivatives of the Eigensurfaces. By considering the local extrema of these surfaces we can locate the largest differences in the corresponding surfaces.

Algorithm 1 delineates these steps and figure 16 shows a sample of the multiresolution process. We consider all of these genes as possible candidates for the cause of the differences in the Eigensurfaces. We then verify these genes using DAVID (the Database for Annotation, Visualization and Integrated Discovery[52]) to check their functionality and choose a gene whose functionality is related to the response of temperature or light stimulus. Figure 17 displays the time series of two such genes. Due to the large amount of data we ran our program on high performance computing facilities of the Keeneland project [53, 54]. Running time of the algorithm on our computing facilities was over 24 hours, where this amount reduced to 3 hours using the Keeneland project facilities.

Algorithms 1.

```
1- A ← interpolated LDHC; B ← interpolated LLHC.
2- Sort A according to similarity of signals; Rearrange B in the same order.
3- eigA ← Eigensurface of A; eigB ← Eigensurface of B (window size 40).
4- D1A ← first derivative of eigA; D1B ← first derivative of eigB.
5- D2A ← second derivative of eigA; D2B ← second derivative of eigB.
6- Delta ← D2A-D2B.
7- E ← The local extrema of Delta.
8- for each point "e" in E, do the following:
    8.1- W2A, W2B ← 40x40 window from A and B that starts from coordinates of e.
    8.2- Delta2 ← difference of second derivatives of Eigensurface of W2A and W2B
            (with sliding window of size 20).
    8.3- E2 ← The local extrema of Delta2.
    8.4- consider "e2" to be maximum of E2.
    8.5- W3A, W3B ← 20x20 window from A and B that starts from coordinates of e2.
    8.6- Delta3 ← difference of second derivatives of Eigensurface of W3A and W3B
            (with sliding window of size 10).
    8.7- E3 ← The local extrema of Delta3.
    8.8- consider "e3" to be maximum of E3.
    8.9- W3A, W3B ← 10x10 window from A and B that starts from coordinates of e3.
    8.10- Delta3 ← difference of second derivatives of Eigensurface of W3A and W3B
            (with sliding window of size 5).
    8.11- E3 ← The local extrema of Delta3.
    8.12- consider "e3" to be maximum of E3.
    8.13- select genes in the 5x5 window that starts from coordinates of e3, as
            possible candidates.
```

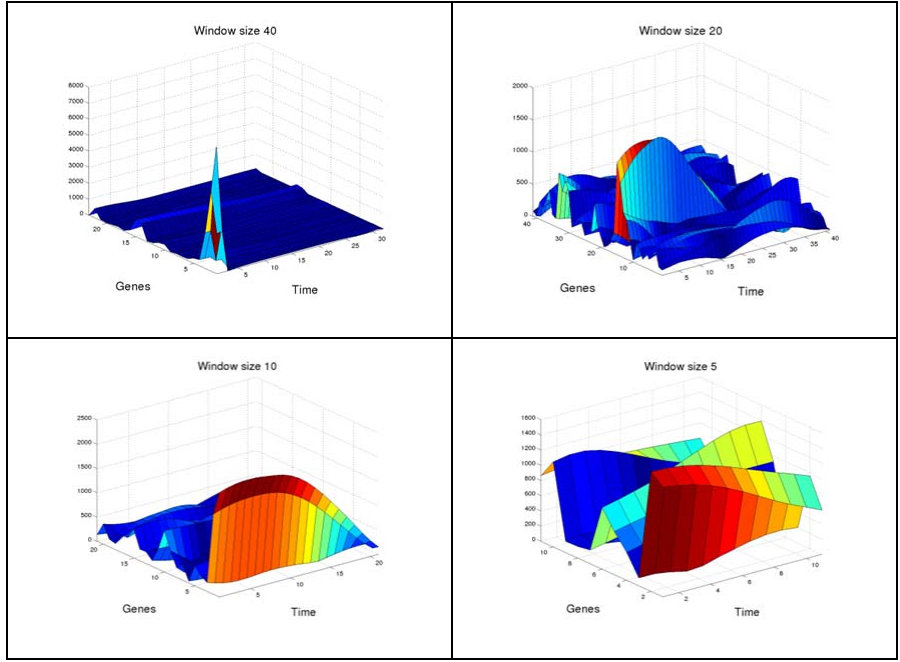

**Fig.16.** This figure shows the multiresolution process. In this process, we first find a local maxima in the Distance surface of Eigensurfaces of windows size 40 (upper left image), and based on that local maxima we choose a small window in the Distance surface of Eigensurfaces of window size 20 (upper right image). By repeating this process using smaller window sizes (bottom left and bottom right), from a 40x40 window in first step, we can reach a 5x5 window in the last step. Therefore we will have smaller number of genes to consider.

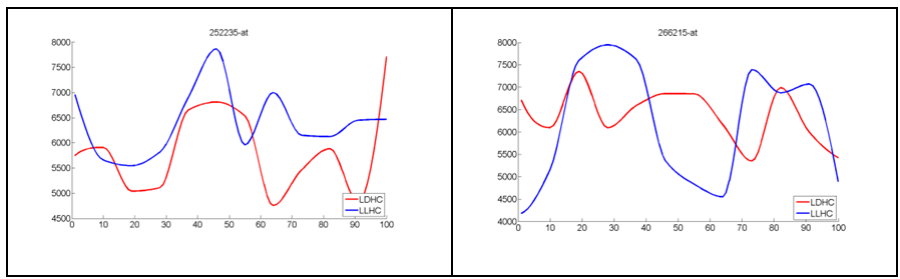

Fig.17. David listed "response to temperature stimulus" and "response to cold" as one of the functionalities of gene AT3G49910 (252235_at), and listed response to "light stimulus" and "response to light intensity" for AT2G06850 (266215_at). The figures show the difference in shape of time series of these two genes in the two different situations. Output of analyzing these data is shown in a supplementary data at (http://www.math.wisc.edu/~dashti/InfoSurf/InfoSurf_Supp1.xls).

**Acknowledgments.** The authors thank Professor Joanne Chory for providing the data sets and discussion about the biological problem. We thank personnel of the "Keeneland: National Institute for Experimental Computing" for their kind supports. This material is based upon work supported by the National Science Foundation under Grant No. 0923296. This projected is partially supported by the National Institute of Health under Grant No. EY21357.